\documentclass{article}
\usepackage{amsmath, amsfonts}
\usepackage{amsmath}
\usepackage{amsfonts}
\usepackage{amssymb}
\usepackage{graphicx}

\begin{document}

\begin{center}
Nonstandard Analysis in Topology:\ \ Nonstandard and Standard Compactifications\ \\[0pt]

by\ \\[0pt]

S. Salbany and T. Todorov\ \\[0pt]\ \\[0pt]

Dedicated to Horst Herrlich on the occasion of his sixtieth
anniversary.\ \\[0pt]\ \\[0pt]
\end{center}

\textbf{Abstract}\ \ Let $\left(  X,T\right)  $ be a topological space, and
$^{*}X$ a non--standard extension of $X$.\ \ There is a natural ``standard''
topology $^{S}T$ on $^{*}X$ generated by $^{*}G$, where $G\in T$.\ \ The
topological space $\left(  ^{*}X,^{S}T\right)  $ will be used to study, in a
systematic way, compactifications of $\left(  X,T\right)  .$\ \newline 



\begin{itemize}
\item [\textbf{1.}]\textbf{Introduction}

Let $\left(  X,T\right)  $ be a topological space and $\left(  ^{\ast}%
X,^{S}T\right)  $ a non--standard enlargement, where $^{S}T$ is generated by
$^{\ast}G$, where $G\in T$.\ \ The space $\left(  ^{\ast}X,^{S}T\right)  $ has
many interesting topological features, but, as may be expected, is very poorly
separated as far as points are concerned.\ \newline 

We show that a wide class of compactifications of $\left(  X,T\right)  $ may
be obtained by rendering $\left(  ^{\ast}X,^{S}T\right)  $ ``separated'', thus
illustrating the usefulness and effectiveness, and broad applicability of the
non--standard compactification $\left(  ^{\ast}X,^{S}T\right)  .$\ \newline 

Conceptually, it is not very common to regard a non--standard model $^{\ast}X$
as a topological space, although this has been done:\ \ A.\ Robinson [16] and
[17], W.A.J. Luxemburg [23],\ H. Gonshor [5],\ L. Haddad [6].\ \ The $\ast
$--open sets have always been part of non--standard techniques , but their
role is more often at the level of the application of transfer principles than
as basic open sets of a topological space $\left(  ^{\ast}X,^{S}T\right)  .$\ \newline 

In this paper we shall study the relationship between topological properties
of $\left(  X,T\right)  $ and their counterparts in $\left(  ^{\ast}%
X,^{S}T\right)  $.\ \ This has led to a unification and, perhaps,
simplification of the exposition concerning compactifications ([22], [14],
[15], [23], [5], [6], [8]).\ \newline 

Arising from considerations related to the Theory of Frames, as well as from
an interest in compactifications that are relevant to theoretical computer
science, there has been an increasing interest in $T_{0}$--compactifications,
described as ``well compacted'' ([20]) and ``stably compact'' ([21]).\ \ We
shall show that these compactifications can also be obtained from
non--standard compactifications in a canonical way.\ \ The relevant reference
for Frames (and also stable compactness) is P.\ Johnstone's book Stone Spaces [10].\ \newline 

A methodological note is appropriate at this stage, concerning the role of the
axiom of choice -- the axiom is essential in topology, to yield the
Stone-\v{C}ech compactification of Tychonoff spaces [12]; it is also essential
in the non--standard approach by providing non--standard enlargements with
adequate saturation ([8], observation before \L\u{o}s Theorem 4.5, Chapter II).\ \newline 

For the topological notions and constructs, we refer to J.L. Kelley [12], and
L. Gillman and M. Jerison [4].\ \ For the relation between standard and
non--standard methods in topology, we refer to L.\ Haddad [6]; for basic
concepts, methods and further developments we refer to T. Lindstr\o m
[13].\ \ For the notions concerning category theory, in particular
reflections, we refer to [1], as well as [10].\ \newline 

\item[\textbf{2.}] \textbf{Non--standard compactifications}

For any topological space $\left(  X,T\right)  $ there is an enlargement
$^{*}X $ which is saturated in the sense that if $\left\{  \left.  F\subseteq
X\right|  i\in I\right\}  $ is a family with the finite intersection property,
then there is $\alpha$ in $^{*}X$ which is in every $^{*}F_{i},\;i\in I$ (see,
for example, [8], Chapter II, \S8).\ \ Thus, $\left(  ^{*}X,^{S}T\right)  $ is
a compact topological space.\ \newline 

The sets $^{*}G$,\ $G\in T$, constitute a base for $^{S}T$.\ \ For reference,
we note that, for $A,B\subseteq X$, we have:

(i) $^{*}\phi=\phi,\;^{*}X=X,\;$(ii)\ $^{*}\left(  A\cup B\right)  =^{*}%
A\cup^{*}B,\;$(iii) $^{*}\left(  A\cap B\right)  =\;^{*}A\cap^{*}B;\;$(iv)
$^{*}X-^{*}A=^{*}\left(  X-A\right)  .$\ \newline 

The monad of $x$ is $\mu\left(  x\right)  =\bigcap\left\{  \left.
^{*}G\right|  G\in T,x\in G\right\}  $.\ \ More generally, for $\alpha\in
^{*}X$,

$\mu\left(  \alpha\right)  =\bigcap\left\{  \left.  ^{*}G\right|  G\in
T,\alpha\in^{*}G\right\}  ,$ similarly we may define $\mu\left(  A\right)  $,
where $A\subseteq^{*}X$.\ \newline 

We recall A.\ Robinson's celebrated criterion for compactness:\ \ $\left(
X,T\right)  $ is compact if and only if $^{*}X=\bigcup\left\{  \left.
\mu\left(  x\right)  \right|  x\in X\right\}  .$\ \newline 

In general $\left(  ^{*}X,^{S}T\right)  $ is not a Hausdorff Topological space
-- the standard open sets are inadequate to separate the rich assortment of
points in $^{*}X$.\ \ We shall provide and example, so as give an idea of the
topological structure of $\left(  ^{*}X,^{S}T\right)  .$\ \newline 

\begin{itemize}
\item [\textbf{2.1}]\textbf{Example\ \ }Consider $\mathbb{N}=\left\{
1,2,3,\ldots\right\}  $ with the topology of upper semi--continuity, basic
open sets being of the form $G_{n}=\left\{  1,2,\ldots,n\right\}
,\;n\in\mathbb{N}$.\ \ Let $\alpha\in^{*}\mathbb{N}$, then $\alpha\in^{*}%
G_{n}$ if and only if $G_{n}\in\alpha$, hence $\alpha$ must be a principal
ultrafilter since $G_{n}$ is finite.\ \ Thus, every $\alpha$ for which
$\left\{  \left.  H\subseteq\mathbb{N}\right|  \alpha\in^{*}H\right\}  $ is a
free ultrafilter, necessarily has only one $^{S}T$--neighbourhood,
$^{*}\mathbb{N}$ itself.\ \ Thus, no ``non standard'' $\alpha$'s can be
separated by $^{S}T$--open sets.\ \newline 

The following observations are important because they give the functoriality
of the $^{*}$--extension.\ \newline 

\item[\textbf{2.2}] \textbf{Proposition\ \ }Let $\left(  X,T\right)  $ and
$\left(  X^{\prime},T^{\prime}\right)  $ be topological spaces and $\left(
^{*}X,^{S}T\right)  ,\;\left(  ^{*}X^{\prime},^{S}T^{\prime}\right)  $
non--standard compactifications.\ \ Then the function $f:\left(  X,T\right)
\rightarrow\left(  X^{\prime},T^{\prime}\right)  $ is continuous if and only
if $^{*}f:\left(  ^{*}X,^{S}T\right)  \rightarrow\left(  ^{*}X^{\prime}%
,^{S}T^{\prime}\right)  $ is continuous.\ \newline 

\textbf{Proof\ \ }For any $G^{\prime}\in T^{\prime}$, we have that
$^{*}G^{\prime}$ is a basic $^{*}T^{\prime}$ open set.\ \ Now $\left(
^{*}f\right)  ^{\leftarrow}\left[  ^{*}G^{\prime}\right]  =(^{*}%
f)^{\leftarrow}\left[  G^{\prime}\right]  $ and the result follows.\ \newline 

\item[\textbf{2.3}] \textbf{Proposition\ \ }Let $f:\left(  X,T\right)
\rightarrow\left(  X^{\prime},J^{\prime}\right)  ,\;g:\left(  X^{\prime
},T^{\prime}\right)  \rightarrow\left(  X^{\prime\prime},T^{\prime\prime
}\right)  ,$ then $^{*}\left(  g\circ f\right)  =\;^{*}g\circ^{*}f.\;\;Also$
$^{*}\left(  \mathbb{I}_{X}\right)  =\mathbb{I}_{^{*}X}.$\ \newline 

We have indicated why $\left(  ^{*}X,^{S}T\right)  $ is a compact topological
space.\ \ In fact, more is true.\ \ Firstly some topological definitions and
their non--standard description.\ \newline 

\item[\textbf{2.4}] \textbf{Definition\ \ }A topological space $\left(
X,T\right)  ,$ not necessarily Hausdorff, is called \textbf{locally compact}
if for every $x\in X$ and open $V\in T$ which is a neighbourhood of $x$, there
is a compact neighbourhood of $x$, $W$, not necessarily open, such that
$W\subseteq V$.\ \newline 

A simple non--standard description of local compactness follows from A.
Robinson's compactness theorem mentioned above.\ \newline 

\textbf{Non--standard Local compactness:\ \ }For every neighbourhood $V$ of a
given $x\in X$, there is a neighbourhood of $x,$ $W$, such that
\[
^{\ast}W\subseteq\bigcup\left\{  \left.  \mu\left(  x\right)  \right|  x\in
W\right\}  \subseteq^{*}V.
\]

The following notion has been called \textbf{supersoberness} ([3], Chapter
VII, 1.10 Definition).\ \ When applied to a compact space $\left(  X,T\right)
$, because it implies a precise form of compactness which specifies, not only
that ultrafilters have adherences, but that these should be of a special form,
we have taken the liberty of naming it \textbf{supercompactness. \ }In [3],
examples illustrating the usefulness of supersoberness may be found in Chapter VII.\ \newline 

\item[\textbf{2.5}] \textbf{Definition\ \ }Let $\left(  X,T\right)  $ be a
topological space, not necessarily $T_{0}.\;\;\left(  X,T\right)  $ is
\textbf{supercompact} if for every ultrafilter $\mathcal{U}$ on $X$ there is
an essentially unique point $x$ such that the set of cluster points of
$\mathcal{U}$, adh $\mathcal{U}$, is the closure of the point $x$:
\[
\text{adh}\mathcal{U}=\;\underset{T}{cl}X.
\]

That $x$ is essentially unique means:\ \ if $x^{\prime}$ is any other point
with that same property, then $x$ and $x^{\prime}$ have precisely the same
$T$--neighbourhoods, i.e. $\mu\left(  x\right)  =\mu\left(  x^{\prime}\right)
$.\ \newline 

\item[\textbf{2.6}] \textbf{Examples}

\begin{itemize}
\item [1.]Every compact Hausdorff space is locally compact and supercompact.

\item[2.] Every supercompact $T_{1}$ space is compact Hausdorff.

\item[3.] Consider $\mathbb{N}_{\infty}=\mathbb{N}\cup\left\{  \infty\right\}
$, where basic $u$--open sets in $\mathbb{N}_{\infty}$ one of the form
$G_{n}=\left\{  0,1,\ldots,n\right\}  $ (thus, the only $u$--neighbourhood of
the point at infinity is $\mathbb{N}_{\infty}$).\ \ Then $\left(
\mathbb{N}_{\infty},u\right)  $ is a $T_{0}$ locally compact supercompact
space.\ \ Indeed, it is the $T_{0}$--locally compact supercompact
compac\-tification $\beta_{2}\left(  \mathbb{N},u\right)  $ of $\left(
\mathbb{N},u\right)  $, see Proposition 4.4.\ \newline 
\end{itemize}

\item[\textbf{2.7}] \textbf{Theorem\ \ }$\left(  ^{*}X,^{S}T\right)  $ is a
locally compact, supercompact enlargement of $\left(  X,T\right)  $.

\textbf{Proof\ \ }We first establish local compactness, by showing that, for
$G\in T$, we have that $^{*}G\;$is $^{S}T$--compact.\ \ Consider a filter
$\mathcal{F}$ of closed sets $^{*}F_{i},\;i\in I$ such that $^{*}F_{i}\cap
^{*}G\neq\phi$, all $i$ in $I$.\ \ Then $\mathcal{F}=\left\{  \left.
F_{i}\cap G\right|  i\in I\right\}  $ is a family of subsets of $X$ which is
closed under finite intersections.\ \ Let $\mathcal{U}$ be an ultrafilter on
$X$ which contains $\mathcal{F}.\;\;$By saturation, there exists $p\in^{*}X$
such that $p\in\bigcap\left\{  \left.  ^{*}\left(  F_{i}\cap G\right)
\right|  i\in I\right\}  $, i.e. $p\in\left(  \bigcap\left\{  \left.
^{*}F_{i}\right|  i\in I\right\}  \right)  \cap^{*}G.\;\;$Thus, $p$ is a
cluster point of $\left\{  \left.  ^{*}F_{i}\cap^{*}G\right|  i\in I\right\}
$ and belongs to $^{*}G$, as required.\ \ To prove supercompactness, let
$\mathcal{V}$ be an ultrafilter on $^{*}X$.\ \ Note that $adh\mathcal{V}%
=\bigcap\left\{  \left.  \bar{H}\right|  H\subseteq^{*}X,\;H\in\mathcal{V}%
\right\}  ,$ so that $adh\mathcal{V}=\bigcap\left\{  \left.  ^{*}F\right|
F\subseteq X,\;F\;\text{is closed, and }^{*}F\in\mathcal{V}\right\}  $.\ \ It
is readily verified that $p=\left\{  \left.  A\subseteq X\right|  \;^{*}%
A\in\mathcal{V}\right\}  $ is an ultrafilter on $X.\;\;$Hence $p\in^{*}%
X$.\ \ We show that $adh\mathcal{V\;}$is the $^{S}T$--closure of $p,$ thus
exhibiting quite explicitly the special minimal point in the
adherence.\ \ Firstly, $p$ is in the adherence of $\mathcal{V}$, since, given
$G\in T$ with $p\in^{*}G$, we have $^{*}G\in\mathcal{V}$, by definition of
$p$.\ \ Hence $^{*}G$ intersects every closed set in $\mathcal{V}$ showing
that $p\in adh\mathcal{V}.\;\;$Let $\alpha\in adh\mathcal{V}$.\ \ If
$\alpha\;$is not in the $^{S}T$--closure of $p$, then there is $G\in T$ such
that $a\in^{*}G$ and $p\notin^{*}G$.\ \ But then $p\in^{*}\left(  X-G\right)
$, hence $^{*}\left(  X-G\right)  \in\mathcal{V}$, so that $^{*}G\cap
^{*}\left(  X-G\right)  =\phi$, contradicting the fact that $\alpha\in
adh\mathcal{V}$.\ \newline 

\textbf{In Summary} -- Every topological space $\left(  X,T\right)  $ may be
embedded into a supercompact locally compact enlargement $\left(  ^{*}%
X,^{S}T\right)  $ with embedding map $\eta_{X}:\left(  X,T\right)
\rightarrow\left(  ^{*}X,^{S}T\right)  $.\ \ The assignment is functorial and
$\eta$ provides, then, a natural transformation from the identity to $*$.\ \newline 

It is natural, and important, to determine the behaviour of $*$ on spaces that
already compact.\ \ We shall examine two cases:\ \ the classical case, even in
the non--standard sense, when $\left(  X,T\right)  $ is compact Hausdorff and
the case where $\left(  X,T\right)  $ is a locally compact, supercompact
$T_{0}$ space.\ \newline 

Before we do so, we shall examine further some separation properties of
$\left(  ^{*}X,^{S}T\right)  .$\ \newline 
\end{itemize}

\item[\textbf{3.}] \textbf{Non--standard compactifications and separation properties.}

We shall describe conditions under which $\left(  ^{*}X,^{S}T\right)  $ is
normal, or regular, or a $T_{0}$--space, in order to illustrate the nature of
$^{S}T$ on $^{*}X$.\ \newline 

\begin{itemize}
\item [\textbf{3.1}]\textbf{Proposition\ \ }$\left(  ^{*}X,^{S}T\right)  $ is
normal if and only if $\left(  X,T\right)  $ is normal.\ \newline 

\textbf{Proof\ \ }Assume that $\left(  ^{*}X,^{S}T\right)  $ is normal and let
$F_{1},F_{2}$ be disjoint closed sets of $\left(  X,T\right)  $.\ \ Then
$^{*}F_{1}$ and $^{*}F_{2}$ are disjoint closed sets of $\left(  ^{*}%
X,^{S}T\right)  $ so, by assumption,\ they can be included in disjoint $^{S}%
T$--open sets with disjoint closures.\ \ Restricting the open sets to $X$
provides two $T$--open sets $G_{1},G_{2}$ with disjoint $T$--closures
containing $F_{1}$ and $F_{2}$, respectively.\ \ Conversely, assume $\left(
X,T\right)  $ is normal.\ \ Let $A,B$ be disjoint closed sets in $\left(
^{*}X,^{S}T\right)  $.\ \ By assumption, $A=\bigcap\left\{  \left.
^{*}F\right|  F\;\text{closed }^{*}F\supseteq A\right\}  ,\;B=\bigcap\left\{
\left.  ^{*}H\right|  H\;\text{closed,\ }^{*}H\supseteq B\right\}  $.\ \ Now
$^{*}F\cap^{*}H=\phi$ for some $^{*}F\supseteq A,\;^{*}H\supseteq B$,
where\ $F,H$ are closed in $X$, otherwise there is $\alpha$ in $^{*}X$ such
that $\alpha\in^{*}\left(  F\cap H\right)  ,$ for all $^{*}F\supseteq
A$,\ $^{*}H\supseteq B$.\ \ This would mean that $\alpha\in A\cap B$, which is
impossible.\ \ Thus we have $F\cap H=\phi$ so there are two disjoint open sets
$G_{1},G_{2}$ such that $F\subseteq G_{1},\;H\subseteq G_{2}.\;\;$Then
$A\subseteq^{*}F\subseteq^{*}G_{1},\;B\subseteq^{*}H\subseteq^{*}G_{2}$ and
$^{*}G_{1}\cap^{*}G_{2}=\phi$, as required.\ \newline 

\item[\textbf{3.2}] \textbf{Proposition\ \ }$\left(  ^{*}X,^{S}T\right)  $ is
regular if and only if every open set in $\left(  X,T\right)  $ is closed.

\textbf{Proof\ \ }Suppose $\left(  X,T\right)  $ has the stated property and
that $\alpha\in^{*}G$ for some $G\in T$.\ \ Since $G$ is open and closed, so
is $^{*}G$, so $\left(  ^{*}X,^{S}T\right)  $ is regular.\ \ Conversely,
assume $\left(  ^{*}X,^{S}T\right)  $ is regular and let $G\in T$.\ \ If $G$
were not closed, then there is $x\in\underset{T}{cl}G-G$.\ \ For each open
neighbourhood of $x,\;H$, we have $H\cap G\neq\phi$, so the family $\left\{
\left.  H\cap G\right|  H\in\mathcal{N}_{x}\right\}  $ has the finite
intersection property.\ \ By saturation, there is $p\in^{*}X$ such that
$p\in^{*}\left(  H\cap G\right)  $, for all $H\in\mathcal{N}_{x}$.\ \ Thus
$p\in\bigcap\left(  ^{*}H\cap^{*}G\right)  ,$ hence $p\in\bigcap\left\{
\left.  ^{*}H\right|  H\in\mathcal{N}_{x}\right\}  \bigcap^{*}G$.\ \ By
regularity, there is $U\in T$ such that $p\in^{*}U\subseteq\underset{^{S}%
T}{cl}^{*}U\subseteq^{*}G$.\ \ But then $X-\underset{T}{cl}U=\left(
^{*}X-\underset{^{S}T}{cl}\;^{*}U\right)  \cap X\in\mathcal{N}_{x}$, since
$x\in\underset{T}{cl}G-G$.\ \ Letting $H=X-\underset{T}{cl}U$, we have $U\cap
H=\phi$, hence $^{*}H\cap^{*}H=\phi$, which contradicts the fact that
$p\in^{*}U$, and $p\in^{*}H$ (since $p\in^{*}W$ for all $W\in\mathcal{N}_{x}$).\ \newline 

\item[\textbf{3.3}] \textbf{Corollary 1\ \ }$\left(  ^{*}X,^{S}T\right)  $ is
a completely regular space (no $T_{0}$ separation assumed) if and only if
every open set in $T$ is closed.

\textbf{Proof\ \ }If $\left(  ^{*}X,^{S}T\right)  $ is completely regular,
then it is regular, hence $\left(  X,T\right)  $ has the desired
property.\ \ Conversely, assume every open set is closed.\ \ Let $G\in T$ and
$\alpha\in^{*}X$ be such that $\alpha\in^{*}G$.\ \ Because $^{*}G$ is open and
closed, there is a continuous real valued function $f:\left(  ^{*}%
X,^{S}T\right)  \rightarrow\left(  \mathbb{R},m\right)  $, where $m$ denotes
the usual topology, such that $G=f^{\leftarrow}\left[  0\right]
,\;X-G=f^{\leftarrow}\left[  1\right]  ,$ the proof is complete.\ \newline 

\item[\textbf{3.4}] \textbf{Corollary 2}\ \ Let $D$ denote the discrete
topology on $\mathbb{N}.\left(  ^{*}\mathbb{N},^{S}D\right)  $ is not a
$T_{0}$--space.

\textbf{Proof\ \ }If $\left(  ^{*}\mathbb{N},^{S}D\right)  $ were $T_{0},$
then it would be $T_{2}$ since $\left(  ^{*}\mathbb{N},^{S}D\right)  $ is
regular, by above.\ \ In which case, since every bounded continuous real
valued function on $\left(  \mathbb{N},D\right)  $ admits an extension to
$\left(  ^{*}\mathbb{N},^{S}D\right)  $ and $\left(  \mathbb{N},D\right)  $ is
dense in $\left(  ^{*}\mathbb{N},^{S}D\right)  $ (trivially, $\left(
X,T\right)  $ is always dense in $\left(  ^{*}X,^{S}T\right)  $, since
$^{*}G\cap X=G$), it would follow that $\left(  ^{*}\mathbb{N},^{S}D\right)  $
is the Stone--\v{C}ech compactification $\beta\left(  \mathbb{N},D\right)  $
of $\left(  \mathbb{N},D\right)  .\;\;$It is well known that this is
impossible (see A. Robinson [16] page 582; or K.D. Stroyan and W.A.J.
Luxemburg [23], 8.1.6,\ 8.1.7,\ 9.1); for a topologist, perhaps the easiest
way is see this is to note that $\beta\left(  \mathbb{N\times N}\right)
\neq\beta\mathbb{N\times\beta N}$ (see, for example, [4]), whereas
$^{*}\left(  \mathbb{N\times N}\right)  =^{*}\mathbb{N\times}^{*}\mathbb{N}$.\ \newline 

\item[\textbf{3.5}] \textbf{Corollary 3\ \ }There is no topology $T$ on
$\mathbb{N}$ for which $\left(  ^{*}\mathbb{N},^{S}T\right)  $ is a $T_{0}$--space.

\textbf{Proof\ \ }Suppose the contrary, that $\left(  ^{*}\mathbb{N}%
\text{,}^{S}T\right)  $ is a $T_{0}$--space for some topology $T$ on
$\mathbb{N}$.\ \ The identity map i:\ \ $\left(  \mathbb{N},D\right)
\rightarrow\left(  \mathbb{N},T\right)  $ is continuous, hence so is its
non--standard extension $^{*}i:\left(  ^{*}\mathbb{N},^{S}D\right)
\rightarrow\left(  ^{*}\mathbb{N},^{S}T\right)  .\;\;$Since $^{*}i$ is
injective and $\left(  ^{*}\mathbb{N},^{S}T\right)  $ is $T_{0}$, it follows
that $\left(  ^{*}\mathbb{N},^{S}D\right)  $ is $T_{0},$ which we know is impossible.\ \newline 

\item[\textbf{3.6}] \textbf{Proposition\ \ }$\left(  ^{*}X,^{S}T\right)  $ is
a $T_{0}$ space if and only if $X$ is finite.

\textbf{Proof\ \ }If $X$ were infinite and $\left(  ^{*}X,^{S}T\right)  $ a
$T_{0}$--space , then there would be a countable subset $\mathbb{N}$ of $X$
with its relative topology, also denoted by $T$, giving $\left(
^{*}\mathbb{N},^{S}T\right)  \subseteq\left(  ^{*}X,^{S}T\right)  $.\ \ Thus
$\left(  ^{*}\mathbb{N},^{S}T\right)  $ is a $T_{0}$ space, which is impossible.\ \newline 

\item[\textbf{3.7}] \textbf{Corollary\ \ }$\left(  ^{*}X,^{S}T\right)  $ is a
Hausdorff space if and only if $\left(  X,T\right)  $ is a finite discrete space.\ \newline 
\end{itemize}

\item[\textbf{4.}] \textbf{Non--standard compactifications of compact spaces
and standard compactifications}\ \newline 

We shall first discuss briefly the compact Hausdorff case.

\begin{itemize}
\item [\textbf{4.1}]\textbf{Proposition\ \ }Let $\left(  X,T\right)  $ be a
compact Hausdorff space and $\left(  ^{*}X,^{S}T\right)  $ a non--standard
extension of $\left(  X,T\right)  .\;\;$There is a continuous retraction
$r_{X}:\left(  ^{*}X,^{S}T\right)  \rightarrow\left(  X,T\right)  ,$ with
$r_{X}$ being the identity when restricted to $X$.\ \newline 

\textbf{Proof\ \ }By Robinson's characterization, $^{*}X=\bigcup\left\{
\left.  \mu\left(  x\right)  \right|  x\in X\right\}  $.\ \ Thus, given
$\alpha\in^{*}X$, there is $x$ such that $\alpha\in\mu\left(  x\right)
.\;\;$Since $\left(  X,T\right)  $ is a Hausdorff space, if $x\neq x^{\prime}%
$, we have $\mu\left(  x\right)  \cap\mu\left(  x^{\prime}\right)  =\phi$,
hence $\alpha\in\mu\left(  x\right)  $ \textbf{for a unique }$x$.\ \ Define
$\underset{X}{r}\left(  \alpha\right)  $ to be that $x$.\ Clearly,
$\underset{X}{r}\left(  x\right)  =x\;$for all $x$ in $X$.\ \ If $G\in T$ and
$x\in G$, then $\alpha\in\mu\left(  x\right)  $ gives $\alpha\in^{*}H$ for all
$H\in T$ that contain $x$.\ \ In particular $\alpha\in^{*}H,$ where $H\in T$
is such that $x\in H\subseteq\bar{H}\subseteq G,\;$($H$ exists by regularity
of $\left(  X,T\right)  $)$.\;\;$If $\beta\in^{*}H$, and $\underset{X}%
{r}\left(  \beta\right)  =x^{\prime},$ then $x^{\prime}\in\bar{H}$, otherwise
$X-\bar{H}$ is an open set containing $x^{\prime}$, hence, by definition of
$\underset{X}{r},\;^{*}\left(  X-\bar{H}\right)  =^{*}X-^{*}\left(  \bar
{H}\right)  $ contains $\beta$, which is impossible since $\beta\in^{*}%
H$.\ \ Thus $\underset{X}{r}$ is a continuous retraction, as stated.\ \newline 

We now show that the Stone--\v{C}ech compactification of a Tychonoff space
$\left(  X,T\right)  $ is simply $\left(  ^{\mathbf{*}}\mathbf{X,}%
^{\mathbf{S}}\mathbf{T}\right)  $ \textbf{made Hausdorff.\ \ }More precisely,
let $\left[  \left(  X,T\right)  \right]  $ denote the $T_{2}$--reflection of
$\left(  X,T\right)  $ ([1]) then $\beta X=\left[  \left(  ^{*}X,^{S}T\right)
\right]  .$\ \newline 

\item[\textbf{4.2}] \textbf{Theorem\ \ }Let $\left(  X,T\right)  $ be a
Tychonoff space.\ Then $\beta\left(  X,T\right)  =\left[  \left(  ^{*}%
X,^{S}T\right)  \right]  .$\ \newline 

\textbf{Proof\ \ }$\left(  X,T\right)  $ is dense in $\left(  ^{*}%
X,^{S}T\right)  $ and the reflection map $\varphi_{X}:\left(  ^{*}%
X,^{S}T\right)  \rightarrow\left[  \left(  ^{*}X,^{S}T\right)  \right]  $ is
continuous, so $\varphi_{X}\left(  X\right)  $ is dense in the compact
Hausdorff space $\left[  \left(  ^{*}X,^{S}T\right)  \right]  .\;\;$Consider
$f:\left(  X,T\right)  \rightarrow\left(  K,S\right)  $ where $\left(
K,S\right)  $ is a compact Hausdorff space.\ We show that there is a map
$F,\;$necessarily unique, such that $F:\left[  \left(  ^{*}X,^{S}T\right)
\right]  \rightarrow\left(  K,S\right)  $ and $F\circ\varphi_{X}\circ\eta
_{X}=f$, where $\eta_{X}:\left(  X,T\right)  \rightarrow\left(  ^{*}%
X,^{S}T\right)  $ is the embedding map of $\left(  X,T\right)  $ into the
non--standard compactification $\left(  ^{*}X,^{S}T\right)  $.\ \ The result
then follows.

Naturality, in the categorical sense, of the constructs is best expressed as a
commutative diagram, given below, from which one can read off the required
$F$.\ \ For convenience, we write $\left[  X\right]  $ in place of $\left[
\left(  X,T\right)  \right]  ,$ and $\left[  f\right]  $ for the reflected
map.\ \newline \ \newline 

$%
\begin{tabular}
[c]{ccccccccccccccc}%
&  &  &  &  &  &  & $\eta_{X}$ &  &  &  &  &  &  & \\
$\left(  X,T\right)  $ &  &  &  &  &  &  &  &  &  &  &  &  &  & $\left(
^{*}X,^{S}T\right)  $\\
&  & $\varphi_{X}$ &  &  &  &  &  &  &  &  &  & $\varphi_{^{*}X}$ &  & \\
&  &  &  &  &  &  & $\left[  \eta_{X}\right]  $ &  &  &  &  &  &  & \\
&  &  &  & $\left[  X\right]  $ &  &  &  &  &  & $\left[  ^{*}X\right]  $ &  &
&  & \\
\multicolumn{1}{l}{} &  &  &  &  &  &  &  &  &  &  &  &  &  & \\
\multicolumn{1}{l}{$f\;\;\;\;\;\;\;$} &  &  &  & \multicolumn{1}{l}{$\left[
f\right]  \;\;\;\;\;$} & \multicolumn{1}{l}{} & \multicolumn{1}{l}{} &
\multicolumn{1}{l}{} & \multicolumn{1}{l}{} & \multicolumn{1}{l}{} &
\multicolumn{1}{r}{$\;\;\;\;\;\left[  ^{*}f\right]  $} &  &  &  &
\multicolumn{1}{r}{$\;\;\;\;\;\;\;^{*}f$}\\
&  &  &  &  &  &  & $\left[  r_{K}\right]  $ &  &  &  &  &  &  & \\
&  &  &  & $\left[  K\right]  $ &  &  &  &  &  & $\left[  ^{*}K\right]  $ &  &
&  & \\
&  &  &  &  &  &  & $\left[  \eta_{K}\right]  $ &  &  &  &  &  &  & \\
&  & $\varphi_{K}$ &  &  &  &  &  &  &  &  &  &  & $\varphi_{^{*}K}$ & \\
&  &  & $\varphi_{K}^{\leftarrow}$ &  &  &  &  &  &  &  &  &  &  & \\
&  &  &  &  &  &  & $r_{K}$ &  &  &  &  &  &  & \\
$\left(  K,S\right)  $ &  &  &  &  &  &  &  &  &  &  &  &  &  & $\left(
^{*}K,^{*}S\right)  $\\
&  &  &  &  &  &  & $\eta_{k}$ &  &  &  &  &  &  &
\end{tabular}
$\ \newline \ \newline 

Observe that $\varphi_{K}$ has an inverse $\varphi_{K}^{\leftarrow},$ since
$\left(  K,S\right)  $ is already Hausdorff.\ \ The required map
is:\ \ $F=\varphi_{K}^{\leftarrow}\circ\left[  r_{K}\right]  \circ\left[
^{*}f\right]  $, since
\begin{align*}
F\circ\varphi_{X}\circ\eta_{X}  & =\varphi_{K}^{\leftarrow}\circ\left[
r_{K}\right]  \circ\left(  \left[  ^{*}f\right]  \cdot\varphi_{X}\right)
\cdot\eta_{X}\\
& =\varphi_{K}^{\leftarrow}\circ\left(  \left[  r_{K}\right]  \cdot
\varphi_{^{*}K}\right)  \circ\left(  ^{*}f\circ\eta_{X}\right) \\
& =\varphi_{K}^{\leftarrow}\circ\varphi_{K}\circ r_{K}\circ\eta_{K}\circ f\\
& =\mathbb{I}_{K}\circ\mathbb{I}_{K}\circ f=f.
\end{align*}
\ \newline 

We now consider $T_{0}$ locally compact supercompact spaces.\ \newline 

\item[\textbf{4.3}] \textbf{Proposition\ \ }Let $\left(  X,T\right)  $ be a
$T_{0}$ locally compact supercompact space and $\left(  ^{*}X,^{S}T\right)  $
a non--standard compactification.\ \ There is a continuous retraction
$r_{X}:\left(  ^{*}X,^{S}T\right)  \rightarrow\left(  X,T\right)  $, which is
the identity when restricted to $X$.\ \newline 

\textbf{Proof}\ \ Let $\alpha\in^{*}X$.\ \ $\alpha$ determines an ultrafilter
$\mathcal{U}_{\alpha}$ on $X$.\ \ As usual, $\mathcal{U}_{\alpha}=\left\{
\left.  A\subseteq X\right|  \alpha\in^{*}A\right\}  $.\ \ By
supercompactness, there is $x$ such that $adh\mathcal{U}_{\alpha}=\underset
{T}{cl}x.\;\;$Define $r_{X}\left(  \alpha\right)  $ to be that $x$, which is,
in fact, unique, as $\left(  X,T\right)  $ is a $T_{0}$ topological
space.\ \ To prove continuity of $r_{X}:\left(  ^{*}X,^{S}T\right)
\rightarrow\left(  X,T\right)  ,$ consider $\alpha\in^{*}X$ and $x=r_{X}%
\left(  \alpha\right)  $.\ \ Given $V\in T$ with $x\in V$, local compactness
ensures that there is $W\in T$ and $K$ compact such that $x\in W\subseteq
K\subseteq V$.\ \ Now $\alpha\in^{*}W$, otherwise $\alpha\in^{*}\left(
X-W\right)  $, so that $X-W\in U_{\alpha}$, contradicting $x\in adh\mathcal{U}%
_{\alpha}$.\ \ Consider now $\beta\in^{*}W$.\ \ We have, then, that $\beta
\in\mu\left(  x^{\prime\prime}\right)  $, for some $x^{\prime\prime}$ in $K$,
since $^{*}W\subseteq\cup\left\{  \left.  \mu\left(  x^{\prime\prime}\right)
\right|  x^{\prime\prime}\in K\right\}  \subseteq^{*}V$.\ \ Hence
$x^{\prime\prime}\in adh\mathcal{U}_{\beta}$, so that $x^{\prime\prime}\in
cl_{T}x^{\prime}$, where $x^{\prime}=r_{X}\left(  \beta\right)  .\;\;$Now
$x^{\prime\prime}\in V$, since $\mu\left(  x^{\prime\prime}\right)
\subseteq^{*}V,$ hence $x^{\prime}\in G$, since $x^{\prime\prime}\in
cl_{T}x^{\prime},$ proving that $r_{X}:\left(  X,T\right)  \rightarrow\left(
^{*}X,^{S}T\right)  $ is continuous.\ \newline 

As in the compact Hausdorff case, an entirely analogous proof will show that
the $T_{0}$ locally compact supercompact reflection of $\left(  X,T\right)
,\;\beta_{2}\left(  X,T\right)  ,$ is ``$\left(  ^{*}\mathbf{X,}^{S}%
\mathbf{T}\right)  $\textbf{\ made }$\mathbf{T}_{0}$''.\ \ More precisely, let
$\left[  \;\;\right]  _{0}$ denote the $T_{0}$--reflector, we have:\ \newline 

\item[\textbf{4.4}] \textbf{Proposition\ \ }Let $\left(  X,T\right)  $ be a
$T_{0}$--space.\ \ Then $\beta_{2}\left(  X,T\right)  =\left[  \left(
^{*}X,^{S}T\right)  \right]  _{0}.$\ \newline 

There is a note of warning that should be mentioned here -- The notion of
``reflection'' that is relevant is the notion of weak reflection of H.
Herrlich, and does not require the uniqueness of the map that solves the
extension problem:\ \ A subcategory $\mathbf{R}$ of a category $\mathbf{A}$ is
reflective, if for every $A$\ in $\mathbf{A}$, there is $R=R\left(  A\right)
$ in $\mathbf{R}$ and $\eta_{A}:A\rightarrow R\left(  A\right)  $ such that if
$A^{\prime}\in\mathbf{A}$ and $f:A\rightarrow A^{\prime}$, then there is $F,$
not necessarily unique, such that $F:R\left(  A\right)  \rightarrow A^{\prime
}$ and
\[
F\circ\eta_{A}=f.
\]
\end{itemize}

\item[\textbf{5.}] \textbf{All Compactifications}

Firstly, a brief reference to compact Hausdorff compactifications of a given
Tychonoff space $\left(  X,T\right)  $.\ \ These may all be obtained by a
uniform method, as described in ([8], page 158).\ \ The method given above
does not refer to continuous real valued functions on $\left(  X,T\right)
$.\ \ However the $T_{2}$ reflection of $\left(  ^{*}X,^{S}T\right)  $ can be
seen to be induced by the $^{*}f^{\prime}s$ where $f$ ranges through the
bounded continuous real valued function, thus establishing a relationship
between the two approaches.\ \newline 

Obtaining all $T_{0}$--locally compact supercompactifications can also be
achieved by considering an analogue of the $Q$--compactifications referred to
above ([8], page 158) -- one considers families of continuous real valued
functions into the Sierpi\'{n}ski dyad $\mathbb{D}=\left\{  0,1\right\}
,\;$with Topology $u=\left\{  \phi,\left\{  0\right\}  ,\left\{  0,1\right\}
\right\}  .$\ \newline 

Lest the reader become too optimistic, it should be mentioned that it is not
possible to obtain \textbf{all} $T_{0}$ compactifications of a $T_{0}$ space
as a type of quotient of $\left(  ^{*}X,^{S}T\right)  $ -- if it were
possible, then the category of compact $T_{0}$ spaces would be (weakly)
reflective in the category of $T_{0}$ spaces, which it is \textbf{not,} as
shown by Miroslav Hu\v{s}ek ([9]; see also [7] for further developments) in
response to a problem posed by Horst Herrlich.\ \newline 
\end{itemize}

\textbf{Acknowledgement:\ }The authors would like to thank the Scuola
Internazional Superiore di Studi Avanzati, Trieste, for its hospitality and
financial support when this work was started.\ \newline The second author
would like to thank Kent Morrison, Ed Glassco and Robert Wolf for a critical
reading of the manuscript and for the correction of many errors.\ \newline 

The authors also wish to thank the referee for his comments, which have
brought about a substantial shortening of the original paper thus allowing it
to focus more clearly on its original intention -- the effectiveness of
non--standard methods when dealing with compactifications.\newpage

\textbf{References}

\begin{itemize}
\item [\lbrack1]] J. Ad\'{a}mek, H. Herrlich and G.E. Strecker,
\textit{Abstract and Concrete Categories}, John Wiley \& Sons, Inc., 1990.

\item[\lbrack2] ] C.C. Chang and H.J Keisler, \textit{Model Theory},
North--Holland, Amsterdam, 1973.

\item[\lbrack3] ] G. Gierz, K.--H. Hofmann, K. Keimel, J.D. Lawson, M.
Mislove, D.S. Scott, \textit{A Compendium of Continuous Lattices},
Springer--Verlag, 1980.

\item[\lbrack4] ] L. Gillman and M. Jerison, \textit{Rings of Continuous
Functions}, Van Nostrand, 1960.

\item[\lbrack5] ] H. Gonshor, \textit{Enlargements containing various kinds of
completions.}\ In \textit{Victoria Symposium on Nonstandard Analysis}, (A.E.
Hurd and P. A. Loeb, editors), Lecture Notes in Mathematics, Vol. 369,
Springer, 1974.

\item[\lbrack6] ] L. Haddad, \textit{Comments on Nonstandard Topology}, Ann.
Sci. Univ. Clermont. S\'{e}r. Math., Fasc. 16 (1978), 1-25.

\item[\lbrack7] ] H. Herrlich, \textit{Compact }$T_{0}$\textit{--spaces and
}$T_{0}$\textit{--compactifications}, Appl. Categ. Structures 1 (1993), 1, 111-132.

\item[\lbrack8] ] A.E. Hurd and P.A. Loeb, \textit{An Introduction to
Nonstandard Real Analysis}, Academic Press, New York, 1985.

\item[\lbrack9] ] M. Hu\v{s}ek, \textit{\v{C}ech--Stone-like compactifications
of general Topological spaces, }Comment.\ \ Math. Univ. Carolinae., 33, 1992, 159--163.

\item[\lbrack10] ] P. Johnstone, \textit{Stone Spaces}, Cambridge Studies in
Advanced Mathematics, Vol 3, CUP, 1982.

\item[\lbrack11] ] H.J. Keisler, \textit{Foundations of Infinitesimal
Calculus}, Prindle, Weber \& Schmidt, 1976.

\item[\lbrack12] ] J.L. Kelley, ``General Topology'', Van Nostrand, 1955.

\item[\lbrack13] ] T. Lindstr\o m, \textit{An Invitation to Nonstandard
Analysis}, In \textit{Nonstandard Analysis and its Applications}, N. Cutland
(Ed), London Mathematical Society Student Texts 10, Cambridge University
Press, Cambridge, 1988, 1--105.

\item[\lbrack14] ] W.A.J. Luxemburg, \textit{A General Theory of Monads}, In
\textit{Applications of Model Theory to Algebra, Analysis, and Probability},
(W.A.J. Luxemburg, ed., 1967). Holt, Rinehard, and Winston, 1969.

\item[\lbrack15] ] M. Machover and J. Hirschfeld, \textit{Lectures on
Non--Standard Analysis}, Lecture Notes in Mathematics, Vol. 94,
Springer--Verlag, 1969.

\item[\lbrack16] ] A. Robinson, \textit{Nonstandard Analysis}, Norht--Holland, 1966.

\item[\lbrack17] ] A Robinson, \textit{Compactification of Groups and Rings
and Nonstandard Analysis}, The Journal of Symbolic Logic, Volume 34, Number 4,
1969, 576--588.

\item[\lbrack18] ] S. Salbany and T. Todorov, \textit{Nonstandard and Standard
Compactifications of Ordered Topological Spaces}, Topology and its
Applications, Vol 47 (1992), 35--52.

\item[\lbrack19] ] S. Salbany, T. Todorov, \textit{Monads and realcompactness}%
, Topology and its Applications, Vol 56 (1994), 99-104.

\item[\lbrack20] ] H. Simmons, \textit{A Couple of Triples}, Top. and Appl.
13, 1982, 201--223.

\item[\lbrack21] ] M.B. Smyth, \textit{Stable Compactifications I}, J. London
Math. Soc. (2) 45 (1992) 321--340.

\item[\lbrack22] ] K.D. Stroyan, \textit{Additional Remarks on the Theory of
Monads}, In \textit{Contributions to nonstandard analysis} (W.A.J. Luxemburg
and A. Robinson, eds), North--Holland, 1972.

\item[\lbrack23] ] K.D. Stroyan and W.A.J. Luxemburg, \textit{Introduction to
the Theory of Infinitesimals}, Academic Press, 1976.
\end{itemize}
\end{document}